\documentclass[11pt]{amsart}

\usepackage{epigamath}

\usepackage[english]{babel}

\input xy
\xyoption {all}

\newtheorem*{theorem}{Theorem}
\newtheorem*{lemma}{Lemma}

\newtheorem*{proposition}{Proposition}

\theoremstyle{definition}
 
\theoremstyle{definition}
\newtheorem{remark}{Remark}

\EpigaVolumeYear{4}{2020} \EpigaArticleNr{3} \ReceivedOn{May 23,
2019}
\InFinalFormOn{January 8, 2020}
\AcceptedOn{February 6, 2020}

\title{On the group of zero-cycles of holomorphic symplectic varieties}
\titlemark{On the group of zero-cycles of holomorphic symplectic varieties}

\author{Alina Marian}
\address{Department of Mathematics, Northeastern University}
\email{a.marian@neu.edu}
\author{Xiaolei Zhao}
\address{Department of Mathematics, University of California, Santa Barbara}
\email{xlzhao@ucsb.edu}

\authormark{A. Marian and X. Zhao}

\AbstractInEnglish{For a moduli space of Bridgeland-stable objects on a K3 surface, we show that the Chow class of a point is determined by the Chern class of the corresponding object on the surface. This establishes a conjecture of Junliang Shen, Qizheng Yin, and the second author.}



\TitleInFrench{Sur le groupe des z\'ero-cycles des vari\'et\'es symplectiques holomorphes}

\AbstractInFrench{Pour un espace de modules d'objets stables au sens de Bridgeland sur une surface K3, nous montrons que la classe de Chow d'un point est d\'etermin\'ee par la classe de Chern de l'objet correspondant sur la surface. Ceci \'etablit une conjecture de Junliang Shen, Qizheng Yi et du second auteur.}

\begin{document}



\maketitle

\begin{prelims}

\DisplayAbstractInEnglish

\bigskip

\DisplayKeyWords

\medskip

\DisplayMSCclass

\bigskip

\languagesection{Fran\c{c}ais}

\bigskip

\DisplayTitleInFrench

\medskip

\DisplayAbstractInFrench

\end{prelims}


\newpage




A Chow-theoretic study of holomorphic symplectic varieties of K3 type was undertaken in \cite{SYZ}. The article considered moduli spaces of Bridgeland-stable sheaf complexes on   a smooth projective K3 surface, and proposed to understand rational equivalence of zero-cycles on the moduli space in terms of the cycle structure of the underlying surface. We prove the following result conjectured in \cite{SYZ}.

\begin{theorem}
 Let $X$ be a smooth projective K3 surface. For a primitive $v \in H^{\star} (X, \mathbb Z)$, and a $v$-generic stability condition $\sigma$, let $\mathsf M_\sigma (v)$ be the moduli space of $\sigma$-stable complexes on $X$ of Mukai vector $v$. Let $F_1$ and $F_2$ be two points in $\mathsf M_\sigma (v)$. Then $$[F_1] = [F_2] \in CH_0 (\mathsf M_\sigma (v) )\iff c_2 (F_1) = c_2 (F_2) \in CH_0 (X). $$
\end{theorem}

\vskip.1in

As noted in \cite{SYZ}, if $F_1$ and $F_2$ have the same Chow class on $\mathsf M_\sigma (v),$ then by restricting a quasi-universal family on $\mathsf M_{\sigma} (v) \times X$ to the two points one concludes that the Chern characters, hence the second Chern classes of $F_1$ and $F_2$ are equal in the Chow ring of $X$. 

We show now that the existence of a quasi-universal family also establishes the other direction. Assume in fact for simplicity that $\mathsf M_\sigma (v)$  admits a universal family $$\mathcal E \to \mathsf M_\sigma (v) \times X.$$ This restriction is certainly not essential, as we will discuss shortly. Denote the two projections by $$\pi: \mathsf M_\sigma (v) \times X \to \mathsf M_\sigma (v) \, \, \, \text{and} \, \, \, \rho: \mathsf M_\sigma (v) \times X \to X.$$ 
The main role in the argument is played by the $\pi$-relative extension complex $$\mathsf W_F = \mathcal Ext_{\pi}^{\bullet}  (\mathcal E, \rho^{\star} F) [1] \, \, \,  \text{on} \, \, \, \mathsf M_\sigma (v),$$
which is well-suited to detect in its Chern classes the Chow class of the point $F \in \mathsf M_\sigma (v).$ For any $E \in \mathsf M_\sigma (v)$ not isomorphic to $F$, by stability and duality we have $$\text{Ext}^0 (E, F) = 0  \, \, \text{and} \,  \, \text{Ext}^2 (E, F) = 0 \,\, \text{on} \, \, X.$$ In this situation, the middle group has constant dimension 
$$\dim \text{Ext}^1 (E, F) = \langle v, \, v\rangle,$$
where the last quantity is the Mukai self-pairing of $v$. Note that the dimension of the moduli space is
$$m = \langle v, \, v \rangle + 2 = \, \dim \, \mathsf M_\sigma (v).$$ On $\mathsf M_\sigma (v) \setminus \{F \}$, we see therefore that $\mathsf W_F$ is a locally free sheaf of rank $m-2$. This implies that on $M_\sigma (v)$ we have: 
\begin{eqnarray*}
c_{m-1} (\mathsf W_F) &=& 0,\\
c_m (\mathsf W_F) &=& \, \text{const} \cdot [F].
\end{eqnarray*}
The constant can be determined by a Thom-Porteous calculation using a locally free complex isomorphic to $\mathsf W_F$. We simply have, in fact: 
\begin{proposition} The top Chern class of the complex $\mathsf W_F$ is 
$$c_m \left (\mathsf W_F \right ) = [F].$$
\end{proposition}
\noindent The proposition gives the theorem via Grothendieck-Riemann-Roch. As $$\text{ch} \, \mathsf W_F = - \pi_{\star} \left (\text{ch} \, \mathcal E^{\vee}  \cdot \rho^{\star} \text{ch} \,F\cdot \rho^{\star} \text{td} X \right),$$ we have in Chow, 
$$\text{ch} \, F_1 = \text{ch} \, F_2 \, \,  \text{on} \, \, X \implies \text{ch} \, \mathsf W_{F_1} = \text{ch} \, \mathsf W_{F_2} \, \, \text{on} \, \, \mathsf M_\sigma (v) \implies c_m (\mathsf W_{F_1} ) = c_m (\mathsf W_{F_2}) \, \, \text{on} \, \, \mathsf M_\sigma (v).$$
The Chern characters on $\mathsf M_\sigma  (v)$ require rational coefficients. As the moduli space has trivial first cohomology group, Roitman's theorem \cite{Bloch}, \cite{Roitman} ensures that $CH_0 (\mathsf M_\sigma (v) )$ is torsion-free, hence the equality of the top Chern classes holds over $\mathbb Z.$ 

The proposition is not a new result, since the extension complex $\mathsf W_F$ repeatedly played a role in the development of moduli theory of sheaves on K3 surfaces. 
Allowing $F$ to vary in moduli, one obtains the complex 
$$\mathsf W = \mathcal Ext^{\bullet}_{\pi} (\mathcal E, \mathcal F) [1] \, \, \, \text{on} \, \, \, \mathsf M_\sigma (v) \, \times \, \mathsf M_\sigma (v),$$ where $\pi: \mathsf M_\sigma (v)\,  \times \,\mathsf M_\sigma (v) \times \, X \to \mathsf M_\sigma (v)\,  \times \, \mathsf M_\sigma (v)$ is the projection, and $\mathcal E, \, \mathcal F$ are the universal objects pulled back to the triple product using the first and second factors respectively. This setting was considered in \cite{markman} when $\mathsf M_\sigma (v)$ parametrizes stable sheaves. In this situation, the main theorem in \cite{markman} asserts

\begin{equation}
\label{diag}
c_m (\mathsf W) = \Delta \in CH_m (\mathsf M_\sigma (v) \times \mathsf M_\sigma (v)).
\end{equation}
\vskip.17in

This diagonal formula holds unsurprisingly for any stability condition $\sigma,$ as we now indicate. The first step in \cite{markman} toward establishing the formula in the sheaf case is to write an explicit geometric three-term complex of locally free sheaves which is isomorphic to the extension complex $\mathsf W$. 
It is clear that for any generic stability condition $\sigma,$ $\mathsf W$ is isomorphic to a three-term complex of locally free sheaves in D$^b(\mathsf M_\sigma  (v) \times \mathsf M_\sigma  (v)),$ 
$$A_0\xrightarrow{\alpha} A_1 \xrightarrow{\beta} A_2 ,$$
where $A_1$ is placed in homological degree $0$. This is true on general grounds \cite[Proposition 5.4]{BM} as long as for all points $x = (E, F) \in \mathsf M_\sigma  (v) \times \mathsf M_\sigma  (v)$, the vanishings
$$H_i(\mathsf W  \overset{L}{\otimes} \mathcal O_x)=0 \, \, \, \text{for} \, \, \,  i\neq-1,0,1$$ hold. As we have 
$\mathsf W \overset{L}{\otimes} \mathcal O_x = \text{Ext}^{\bullet}_X (E, F) [1],$
and  $E, \, F$ both belong to the heart of the stability condition $\sigma$, this is true. 

Note next that for the highest extension sheaf $\mathcal Ext^2_{\pi} (\mathcal E, \mathcal F),$ cohomology and base change commute, so this sheaf is a line bundle supported on the diagonal $\Delta$. By contrast, $\mathcal Ext^0_{\pi} (\mathcal E, \mathcal F)$ is torsion-free hence zero. For the above locally free complex, the map $\alpha$ is thus injective, while the cokernel of $\beta$ is a line bundle supported on the diagonal. Furthermore, Grothendieck-Verdier duality gives $$\mathsf W^{\vee}= \mathcal Ext_{\pi}^{\bullet}  (\mathcal F, \mathcal E) [1],$$ so the cokernel sheaf of $\alpha^{\vee}: A_1^{\vee}  \to A_0^{\vee}$ is also supported as a line bundle on $\Delta.$

The following general lemma proven in \cite{markman} by a suitable Thom-Porteous calculation then yields the diagonal formula \eqref{diag} for any generic stability condition $\sigma$.
\begin{lemma}[\cite{markman}, Lemma 4]
Let $$A_0\xrightarrow{\alpha} A_1 \xrightarrow{\beta} A_2 $$ be a complex of locally free sheaves on a smooth projective variety $M$, such that
\begin{itemize}
\item[(1)] The cohomology sheaves satisfy $\mathcal H_0 = 0$ and $\mathcal H_2$ is a line bundle supported on a smooth subvariety $Z$ of even codimension $m$. 
\item[(2)] The cokernel sheaf of $\alpha^{\vee}: A_1^{\vee}  \to A_0^{\vee}$ is also supported as a line bundle on $Z$.
\item[(3)] The ranks satisfy $a_1 - a_0 - a_2 = m-2 \geq 0.$
\end{itemize}
Then $c_m (A_1 - A_0 - A_2) = [Z].$
\label{lemma}
\end{lemma} 

\noindent Since for any $F \in \mathsf M_\sigma  (v)$, the complex $\mathsf W_F$ is the (derived) pullback $$\mathsf W_F = \iota^{\star}_{F} \mathsf W$$ under the inclusion $$\iota_F:  \mathsf M_\sigma (v) \simeq  \mathsf M_\sigma  (v) \times F \hookrightarrow  \mathsf M_\sigma (v) \times  \mathsf M_\sigma (v),$$
the proposition is simply obtained from the diagonal formula \eqref{diag} by pullback of $c_m (\mathsf W).$

\vskip,1in

The moduli space $\mathsf M_\sigma (v)$ does not generally admit a universal family. Nevertheless, as explained in  \cite[Theorem A.5]{mukai}, and later \cite[Section 3] {markman}, a universal object can be glued together over $\mathbb P \times X$, where $\mathbb P \to  \mathsf M_\sigma (v)$ is a suitable projective bundle of a certain rank $N$. As further shown in \cite[Lemma 4.2] {kls}, for any point $F \in \mathsf M_\sigma (v)$, by taking the intersection of $\mathbb P$ with $N$ generic hyperplanes, one obtains 
 a smooth projective variety $\widetilde M$ parametrizing a family $$\mathcal E \to \widetilde M \times X$$ of $\sigma$-stable objects in D$^b(X)$, such that the classifying morphism $f:\widetilde M\to \mathsf M_\sigma  (v)$ is surjective, generically finite, and \'etale over a neighborhood of $F$.

Considering the $\pi$-relative extension complex
$$\mathsf W_F = \mathcal Ext_{\pi}^{\bullet}  (\mathcal E, \rho^{\star} F) [1] \, \, \,  \text{on} \, \, \, \widetilde M,$$
the argument above gives  
$$f_*c_m \left (\mathsf W_F \right ) = n [F],$$ where $n$ is the generic degree of $f$. 
Since $CH_0 (\mathsf M_\sigma (v) )$ is torsion-free, the theorem follows from this formula in conjunction with Grothendieck-Riemann-Roch. \qed

\vskip.3in

\begin{remark}
In the classic paper \cite{mukai}, for isotropic Mukai vectors $v$, the author used $\mathsf W_F$ to establish irreducibility of the two-dimensional moduli space $\mathsf M_\sigma (v).$ An argument with the same flavor was later made in \cite{kls} to show irreducibility in a more general setting. In a different direction, the aim in \cite{markman} is to represent the diagonal in terms of Chern classes of the universal sheaves, which entails that the cohomology of the moduli space is generated by the K\"{u}nneth components of the universal Chern classes. The origin of this idea can in turn be traced back to the earlier work \cite{beauville}, \cite{es}.

The subring of $CH_{\star} \left (\mathsf M_\sigma (v) \right) $ generated in a suitable sense by the universal Chern classes is important for understanding the cycle class map to cohomology along the lines of \cite{beauville1}, \cite{Voisin_conj}. 

\end{remark}

\begin{remark}
The argument above gives the same statement for abelian surfaces. For surfaces with effective anti-canonical divisors (when Ext$^2$ vanishes by Serre duality), the proposition reads
$$c_m \left (\mathsf W_F \right ) = (-1)^m [F],$$
and the proof is completely similar. 
\end{remark}

\begin{remark}
The consequences of the theorem were investigated in \cite{SYZ}. In particular, the result implies the existence of algebraic coisotropic subvarieties $q_i: Z_i \dashrightarrow B_i$ for every $i$, where $Z_i$ is a subvariety of codimension $i$ in $\mathsf M_\sigma (v)$, and the generic fibers of $q_i$ are constant cycle subvarieties of dimension $i$. The existence of such varieties was conjectured in \cite{Voisin} for any holomorphic symplectic variety $M$, and was used to construct a Beauville-Voisin filtration on $CH_0 (M)$. Another filtration on $CH_0 (\mathsf M_\sigma (v))$ was constructed in \cite{SYZ}, and was shown to be related to the one in \cite{Voisin} assuming the theorem.
\end{remark}

\section*{Acknowledgements}

We thank Dragos Oprea for useful related discussions. A. M. was supported by the NSF through grant DMS 1601605. She thanks Davesh Maulik and the MIT mathematics department for their hospitality during the fall semester 2017. X. Z. is partially supported by the Simons Collaborative Grant 636187.


\end{document}